\documentclass{amsart}

\usepackage{amsmath,amssymb}
\usepackage[latin1]{inputenc}
\usepackage[dvips]{graphics}
\input{xy}
\xyoption{poly}
\xyoption{2cell}
\xyoption{all}
\usepackage{epsfig}  
\usepackage{color}

\newcommand{\za}{\alpha}
\newcommand{\zb}{\beta}

\newcommand{\zD}{\Delta}

\newcommand{\zg}{\gamma}

\newcommand{\zs}{\sigma}
\newcommand{\zS}{\Sigma}

\newcommand{\Sbar}{\overline{S}}

\newtheorem{thm}{Theorem}[section]
\newtheorem*{thm*}{Theorem \ref{thm main}}
\newtheorem*{thm**}{Theorem \ref{thm y}}
\newtheorem*{thmF}{Theorem \ref{thm F}}
\newtheorem{prop}[thm]{Proposition}

\newtheorem{cor}[thm]{Corollary}
\newtheorem{lem}[thm]{Lemma}

\newtheorem{definition}{Definition} 
\newtheorem{rem}[thm]{Remark}

\newenvironment{pf}{{Proof}.}

\title{Cluster expansion formulas and perfect matchings}
\author{Gregg Musiker and Ralf Schiffler}
\date{\today}

\begin{document}

\begin{abstract} 
We study cluster algebras with principal coefficient systems that are associated to unpunctured surfaces. We give a direct formula for the Laurent polynomial expansion of cluster variables in these cluster algebras in terms of perfect matchings of a certain graph $G_{T,\gamma}$ that is constructed  from the surface by recursive glueing of elementary pieces that we call tiles. 
We also give a second formula for these Laurent polynomial expansions in terms of subgraphs of the graph $G_{T,\gamma}$.

\end{abstract}

\maketitle
%%%%%%%%%%%%%%%%%%%%%%%%%%%%%%%%
%% 
%%  SECTION
%%
%%%%%%%%%%%%%%%%%%%%%%%%%%%%%%%%

\begin{section}{Introduction}\label{section intro}
Cluster algebras, introduced in \cite{FZ1}, are commutative algebras
equipped with a distinguished set of generators, the \emph{cluster
  variables}. The cluster variables are grouped into sets of constant
cardinality $n$, the \emph{clusters}, and the integer $n$ is called
the \emph{rank} of the cluster algebra. Starting with an initial cluster
$\mathbf{x}$ (together with a skew symmetrizable integer $n\times n$ matrix $B=(b_{ij})$ and a coefficient vector
$\mathbf{y}=(y_i)$ whose entries are elements of a torsion-free
abelian group $\mathbb{P}$) the set of cluster variables is obtained
by repeated application of so called \emph{mutations}. 
To be more precise, let $\mathcal{F}$ be 
 the field of rational functions in the indeterminates $x_1,x_2,\ldots,x_n$
over the quotient field of the integer group ring
$\mathbb{ZP}$. Thus $\mathbf{x}=\{x_1,x_2,\ldots,x_n\}$ is a
transcendence basis for $\mathcal{F}$.
For every $k=1,2,\ldots,n$, the  mutation
$\mu_k(\mathbf{x})$ of the cluster
$\mathbf{x}=\{x_1,x_2,\ldots,x_n\}$ is a new cluster 
$\mu_k(\mathbf{x})=\mathbf{x}\setminus \{x_k\}\cup\{x_k'\}$ obtained
  from $\mathbf{x}$ by replacing the cluster variable $x_k$ by the new
  cluster variable 
\begin{equation}\label{intro 1}
x_k'= \frac{1}{x_k}\,\left(y_k^+\,\prod_{b_{ki}>0} x_i^{b_{ki}} +
y_k^-\,\prod_{b_{ki}<0} x_i^{-b_{ki}}\right)
\end{equation}  
in $\mathcal{F}$, where $y_k^+,y_k^-$ are certain monomials in $y_1,y_2,\ldots,y_n$.
Mutations also change the attached matrix $B$ as well as the coefficient
vector $\mathbf{y}$, see \cite{FZ1}.

The set of all cluster variables is the union of all clusters obtained
from an initial cluster $\mathbf{x}$ by repeated mutations. Note that
this set may be infinite.

It is clear from the construction that every cluster variable is a
rational function in the initial cluster variables
$x_1,x_2,\ldots,x_n$. In \cite{FZ1} it is shown that every cluster
variable $u$ is actually a Laurent polynomial in the $x_i$, that is,
$u$ can be written as a reduced fraction 
\begin{equation}\label{intro 2}
u=\frac{f(x_1,x_2,\ldots,x_n)}{\prod_{i=1}^n x_i^{d_i}},
\end{equation}  
where $f\in\mathbb{ZP}[x_1,x_2,\ldots,x_n]$ and $d_i\ge 0$.
The right hand side of equation (\ref{intro 2}) is called the
\emph{cluster expansion} of $u$ in $\mathbf{x}$.

The cluster algebra is determined by the initial matrix $B$ and the choice of the coefficient system. 
A canonical choice of coefficients is the {\em principal} coefficient system, introduced in \cite{FZ4}, which means that the coefficient group $\mathbb{P}$ is the free abelian group on $n$ generators $y_1,y_2,\ldots,y_n$, and the initial coefficient tuple $\mathbf{y}=\{y_1,y_2,\ldots,y_n\}$ consists of these $n$ generators. In \cite{FZ4}, the authors show that knowing the expansion formulas in the case where the cluster algebra has principal coefficients allows one to compute the expansion formulas for arbitrary coefficient systems.

Inspired by the work of Fock and Goncharov \cite{FG1,FG2,FG3} and
Gekhtman, Shapiro and Vainshtein \cite{GSV1,GSV2} which discovered
cluster structures in the context of
Teichm\"uller theory, Fomin, Shapiro and Thurston  \cite{FST,FT}  initiated a
systematic study of the cluster algebras arising from 
triangulations of a surface with boundary and marked points.
 In this approach,  cluster variables in the cluster algebra correspond to arcs in the surface, and clusters
correspond to triangulations. 
In \cite{S2}, building on earlier results in \cite{S1,ST}, this model was used to give a direct expansion formula
for cluster variables in cluster 
 algebras associated to unpunctured
surfaces, with arbitrary coefficients, in terms of
certain paths on the triangulation. 

Our first main result in this paper is a new parametrization of this formula
in terms of perfect matchings of a certain weighted graph that is constructed
from the surface by recursive glueing of elementary pieces that we
call tiles.
To be more precise, let $x_\gamma$ be a cluster variable corresponding
to an arc $\zg$ in the unpunctured surface and let $d$ be the
number of crossings between $\zg$ and the triangulation $T$ of the
surface. Then $\zg$ runs through $d+1$ triangles of $T$
and each pair of consecutive triangles forms a
quadrilateral which we call a tile. So we obtain $d$ tiles, each of
which is a weighted graph, whose weights are given by the 
cluster variables $x_\tau$ associated to the arcs $\tau$ of the
triangulation $T$. 

 We obtain a weighted graph $G_{T,\zg}$ by 
glueing the $d$ tiles in a specific way and then
deleting the diagonal in each tile.
 To any perfect matching
$M$ of this 
graph we associate its weight $w(M)$ which is  the product of the
weights of its edges, hence a product of cluster variables.
 We prove the following cluster expansion formula:

\begin{thm*} 
\[x_\zg = \sum_M \frac{w(M)\,y(M)}{x_{i_1}x_{i_2}\ldots x_{i_d}}, \]
where the sum is over all perfect matchings $M$ of $G_{T,\zg}$,  $w(M)$ is the weight of $M$, and $y(M)$ is a monomial in $\mathbf{y}$.
\end{thm*}

We also give a formula for the coefficients $y(M)$ in terms of perfect matchings as follows. The $F$-polynomial $F_\zg$, introduced in \cite{FZ4} is 
 obtained from the Laurent polynomial $x_\zg$ (with
principal coefficients) by 
substituting $1$ for each of the cluster variables $x_1,x_2,\ldots,x_n$.
By \cite[Theorem 6.2, Corollary 6.4]{S2}, the $F$-polynomial
has constant term $1$ and a unique term of maximal degree that is divisible by all the other occurring monomials. 
The two corresponding matchings  are the unique two matchings that have all their edges on the boundary of the graph $G_{T,\zg}$. 
We denote by $M_-$ the one with $y(M_-)=1$ and the other by $M_+$.
Now, for an arbitrary perfect matching $M$, the coefficient $y(M)$ is determined by the set of edges of the symmetric difference  
$M_-\ominus M =(M_-\cup M)\setminus (M_-\cap M)$ as follows.

\begin{thm**}
The set $M_-\ominus M$ is the set of boundary edges of a (possibly disconnected) subgraph $G_M$ of $G_{T,\zg}$ which is a union of tiles
$ G_M =\cup_{j\in J} S_j. $
Moreover,
\[y(M)=\prod_{j\in J} y_{i_j}\]
\end{thm**}  

As an immediate corollary, we see that the corresponding $g$-vector, introduced in \cite{FZ4}, is  
\[g_\zg= \deg\left(\frac{w(M_-)}{x_{i_1}\cdots x_{i_d}}\right).\]

Our third main result is yet another description of the formula of Theorem \ref{thm main} in terms of the graph $G_{T,\zg}$ only. In order to state this result, we need some notation. If $H$ is a graph,  let $c(H)$ be the number of connected components of  $H$, 
let $E(H)$ be the set of edges of $H$, and denote by $\partial H$ the set of boundary edges of $H$.
Define $\mathcal{H}_k$ to be 
the set of all subgraphs $H$ of $G_{T,\zg}$ such that $H$ is a union of $k$ tiles $H=S_{j_1}\cup\cdots\cup S_{j_k}$ and such that
the number of edges of $M_-$ that are contained in $H$ is equal to $k+c(H)$.
For $H\in \mathcal{H}_k$, let \[y(H)=\prod_{S_{i_j} \textup{\,tile\,in\,}H} y_{i_j}.\]

\begin{thmF}
The cluster expansion of the cluster variable $x_\zg$ is given by
\[x_\zg=\sum_{k=0}^d \  \sum_{H\in \mathcal{H}_k} \frac{w(\partial H\ominus M_-)\,y(H)}{x_{i_1} x_{i_2}\cdots x_{i_d}}.\]
\end{thmF}

Theorem \ref{thm main} has interesting intersections with work of other
people. In \cite{CCS2}, the authors obtained a formula for the
denominators of the cluster expansion in types $A,D$ and $E$, see also
\cite{BMR}. In \cite{CC,CK,CK2} an expansion formula was given in the
case where the cluster algebra is acyclic and the cluster lies in an
acyclic seed. Palu generalized this formula to arbitrary
clusters in an acyclic cluster algebra \cite{Palu}.  These formulas
use the cluster category introduced in \cite{BMRRT}, and in \cite{CCS1} for
type $A$, and do not give information about the coefficients. 

Recently, Fu and Keller generalized this formula further to cluster algebras  with principal coefficients that admit a categorification by a 2-Calabi-Yau category \cite{FK}, and, combining results of \cite{A} and \cite{ABCP,LF}, such a categorification exists in the case of cluster algebras associated to unpunctured surfaces.

In \cite{SZ,CZ,Z,MP} cluster expansions for cluster algebras of
rank 2 are given, in \cite{Propp,CP,FZ3} the case $A$ is
considered.  
In section 4 of \cite{Propp}, Propp describes two constructions of snake graphs, the latter of which are unweighted analogues for the case A of the graphs $G_{T,\gamma}$ that we present in this paper.  Propp assigns a snake graph to each arc in the triangulation of an $n$-gon and shows that the numbers of matchings in these graphs satisfy the Conway-Coxeter frieze pattern induced by the Ptolemy relations on the $n$-gon.
In \cite{M} a cluster expansion for cluster algebras of
classical 
type is given for clusters that lie in a bipartite seed, 
and the forthcoming work of \cite{MW} will concern cluster expansions for cluster algebras of classical type with principal coefficients, for an arbitary seed.

The formula for $y(M)$ given in  
Theorem \ref{thm y} also can be formulated in terms of height functions, as found in literature such as \cite{EKLP} or \cite{ProppLattice}.  We discuss this connection in Remark \ref{height} of section \ref{sect y}.

The paper is organized as follows. In section \ref{sect FST}, we
recall the construction of cluster algebras from surfaces of \cite{FST}.
Section \ref{sect main} contains the construction of the graph $G_{T,\zg}$ and the statement of the cluster expansion formula. Section \ref{sect proof} is devoted to the proof of the expansion formula.
 The formula for $y(M)$ and the formula for the $g$-vectors is given in section \ref{sect y}.
In section \ref{sect F-polynomial}, we present the expansion formula in terms of subgraphs and deduce a formula for the $F$-polynomials. We give an example in section \ref{sect example}.

\vspace{2em} {\bf Acknowledgements.} The authors would like to thank Jim Propp for useful conversations related to this work.

\end{section} 
%%
%% 
%%

%%%%%%%%%%%%%%%%%%%%%%%%%%%%%%%%%%%%%%%%%%%%%%%%%%%%%%%%%%%%%%%%%%%%%%%%
%%
%%                            SECTION 1
%%
%%%%%%%%%%%%%%%%%%%%%%%%%%%%%%%%%%%%%%%%%%%%%%%%%%%%%%%%%%%%%%%%%%%%%%%%

%%
%% 
%%

\begin{section}{Cluster algebras from surfaces}\label{sect FST}
In this section, we recall the construction of \cite{FST} in the
case of surfaces without punctures.

Let $S$ be a connected oriented 2-dimensional Riemann surface with
boundary and $M$ a non-empty set of marked points in the closure of
$S$ with at least one marked point on each boundary component. The
pair $(S,M)$ is called \emph{bordered surface with marked points}. Marked
points in the interior of $S$ are called \emph{punctures}.  

In this paper we will only consider surfaces $(S,M)$ such that all
marked points lie on the boundary of $S$, and we will refer to $(S,M)$
simply by \emph{unpunctured surface}. 

We say that two curves in $S$ \emph{do not cross} if they do not intersect
each other except that endpoints may coincide.

\begin{definition}
An \emph{arc} $\zg$ in $(S,M)$ is a curve in $S$ such that 
\begin{itemize}
\item[(a)] the endpoints are in $M$,
\item[(b)] $\zg$ does not cross itself,
\item[(c)] the relative interior of $\zg$ is disjoint from $M$ and
  from the boundary of $S$,
\item[(d)] $\zg$ does not cut out a monogon or a digon. 
\end{itemize}   
\end{definition}     
 Curves that connect two
marked points and lie entirely on the boundary of $S$ without passing
through a third marked point are called \emph{boundary arcs}.
Hence an arc is a curve between two marked points, which does not
intersect itself nor the boundary except possibly at its endpoints and
which is not homotopic to a point or a boundary arc.

Each arc is considered up to isotopy inside the class of such curves. 
 Moreover, each arc is considered up to orientation, so if an arc has endpoints $a,b\in M$ then it can be represented by a curve that runs from $a$ to $b$, as well as by a curve that runs from $b$ to $a$.

For any two arcs $\zg,\zg'$ in $S$, let $e(\zg,\zg')$ be the minimal
number of crossings of $\zg$ and $\zg'$, that is, $e(\zg,\zg')$ is the
minimum of
the numbers of crossings of  arcs $\za$ and $\za'$, where $\za$ is
isotopic to $\zg$ and $\za'$ is isotopic to $\zg'$.
Two arcs $\zg,\zg'$ are called \emph{compatible} if $e(\zg,\zg')=0$. 
A \emph{triangulation} is a maximal collection of
compatible arcs together with all boundary arcs. 
The arcs of a 
triangulation cut the surface into \emph{triangles}.
Since $(S,M)$ is an unpunctured surface, the three sides of each
triangle are distinct (in contrast to the case of surfaces with
punctures).  Any triangulation  has
$n+m$ elements, $n$ of which  are arcs in $S$, and the remaining $m$
elements are boundary arcs. Note that the number of boundary arcs
is equal to the number of marked points.

\begin{prop}\label{prop rank}
The number $n$ of arcs in any triangulation is  given by the formula 
$n=6g+3b+m-6$,  where $g$ is the
genus of $S$, $b$ is the number of boundary components and $m=|M|$ is the
number of marked points. The number $n$ is called the \emph{rank} of $(S,M)$.
\end{prop}  
\begin{pf} \cite[2.10]{FST}
\qed
\end{pf}  

Note that $b> 0$ since the set $M$ is not empty.
 Table \ref{table 1} gives some examples of unpunctured surfaces.

\begin{table}
\begin{center}
  \begin{tabular}{ c | c | c || l  }
  \  b\ \  &\ \  g \ \   & \ \  m \ \  &\  surface \\ \hline
    1 & 0 & n+3 & \ polygon \\ 
    1 & 1 & n-3 & \ torus with disk removed \\
    1 & 2 & n-9 & \ genus 2 surface with disk removed \\\hline 
    2 & 0 & n & \ annulus\\
    2 & 1 & n-6 & \ torus with 2 disks removed \\ 
    2 & 2 & n-12 & \ genus 2 surface with 2 disks removed \\ \hline
    3 & 0 & n-3 & \ pair of pants \\ \\
  \end{tabular}
\end{center}
\caption{Examples of unpunctured surfaces}\label{table 1}
\end{table}

Following  \cite{FST}, we associate a cluster algebra
to the unpunctured surface $(S,M)$ as follows.
 Choose any triangulation
$T$, let $\tau_1,\tau_2,\ldots,\tau_n$ be the $n$ interior arcs of
$T$ and  denote the  $m$ boundary
arcs of the surface by $\tau_{n+1},\tau_{n+2},\ldots,\tau_{n+m}$. 
For any triangle $\Delta$ in $T$ define a matrix 
$B^\Delta=(b^\Delta_{ij})_{1\le i\le n, 1\le j\le n}$  by
\[ b_{ij}^\Delta=\left\{
\begin{array}{ll}
1 & \textup{if $\tau_i$ and $\tau_j$ are sides of 
  $\Delta$ with  $\tau_j$ following $\tau_i$  in the }
\\ &\textup{ counter-clockwise order;}\\
-1 &  \textup{if $\tau_i$ and $\tau_j$ are sides of
  $\Delta$ with  $\tau_j$ following $\tau_i$  in the }\\
&\textup{  clockwise  order;}\\
0& \textup{otherwise.}
\end{array} \right. \]
Then define the matrix 
$ B_{T}=(b_{ij})_{1\le i\le n, 1\le j\le n}$  by
$b_{ij}=\sum_\Delta b_{ij}^\Delta$, where the sum is taken over all
triangles in $T$. Note that the boundary arcs of the triangulation are ignored in the definition of $B_{T}$.
Let $\tilde B_{T}=(b_{ij})_{1\le i\le 2n, 1\le j\le n}$ be the
$2n\times n$ matrix whose upper $n\times n$ part is $B_{T}$ and whose
lower $n\times n$ part is the identity matrix.
The matrix $B_{T}$ is skew-symmetric and each of its entries  $b_{ij}$ is either
$0,1,-1,2$, or $-2$, since every arc $\tau$ can be in at most two triangles. 
An example where $b_{ij}=2 $ is given in Figure
\ref{fig bij=2}.

\begin{figure}
\centering
\input{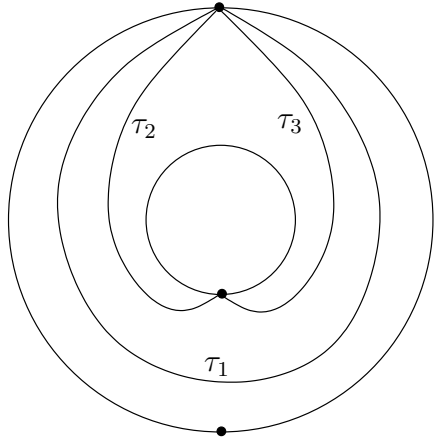}
\caption{A triangulation with $b_{23}=2$ \label{fig bij=2}}
\end{figure}

Let $\mathcal{A}(\mathbf{x}_{T},\mathbf{y}_{T},B_{T})$ be the cluster algebra with principal
coefficients in the triangulation $T$, that is,  $\mathcal{A}(\mathbf{x}_{T},\mathbf{y}_{T},B_{T})$ is
given by the seed
$(\mathbf{x}_{T},\mathbf{y}_{T},B_{T})$ where
$\mathbf{x}_{T}=\{x_{\tau_1},x_{\tau_2},\ldots,x_{\tau_n}\}$ is the
cluster associated to the triangulation $T$, and the initial coefficient
vector $\mathbf{y}_{T}=(y_1,y_2,\ldots,y_n)$ is the vector of generators
of $\mathbb{P}=\textup{Trop}(y_1,y_2,\ldots,y_n)$. 

For the  boundary arcs we define $x_{\tau_k}=1$,  $k=n+1,n+2,\ldots,n+m$.

For each $k=1,2,\ldots,n$,  there is a unique quadrilateral  in $T\setminus \{\tau_k\}$
in which $\tau_k$ is one of the diagonals. Let $\tau_k'$ denote the
other diagonal in that quadrilateral.
Define the
\emph{flip} $\mu_k T$ to be the triangulation
$T\setminus\{\tau_k\}\cup\{\tau_k'\}$.
The mutation $\mu_k$ of the seed $\zS_T$ in the cluster algebra
$\mathcal{A}$ corresponds to the flip $\mu_k$ of the triangulation $T$
in the following sense.
The matrix $\mu_k(B_T)$ is the matrix corresponding to the
triangulation $\mu_k T$,
the cluster
$\mu_k(\mathbf{x}_T)$ is $\mathbf{x}_T\setminus\{x_{\tau_k}\}\cup
\{x_{\tau_k'}\}$, and the corresponding exchange relation is given by
\[x_{\tau_k} x_{\tau_k'} = x_{\rho_1} x_{\rho_2} y^+ +  x_{\zs_1}
x_{\zs_2} y^-,
\]
where $y^+,y^-$ are some coefficients, and $\rho_1,\zs_1,\rho_2,\zs_2$
are the sides of the quadrilateral in which $\tau_k$ and $\tau_k'$ are
the diagonals, such that  $\rho_1,\rho_2$ are opposite sides and
$\zs_1,\zs_2$ are opposite sides too.

\end{section}

%%%%%%%%%%%%%%%%%%%%%%%%%%%%%%%%
%% 
%%  SECTION
%%
%%%%%%%%%%%%%%%%%%%%%%%%%%%%%%%%

\begin{section}{Expansion formula}\label{sect main}

In this section, we will present an expansion formula for the cluster variables in terms of perfect matchings of a graph that is constructed recursively using so-called \emph{tiles}.

\begin{subsection}{Tiles}\label{sect tiles}
For the purpose of this paper, a \emph{tile} $\Sbar_k$ is a planar four vertex
  graph with five weighted edges having the shape of two equilateral
  triangles that share one edge, see Figure \ref{figtile}. 
\begin{figure}
\input{figtile.pstex_t}
\caption{The tile $\Sbar_k$}\label{figtile}
\end{figure}   
The weight on each edge of the tile $\Sbar_k$ is a single variable. The
unique interior edge is called \emph{diagonal} and the four exterior
edges are called \emph{sides} of $\Sbar_k$.
We shall use $S_k$ to denote the graph obtained from $\Sbar_k$ by
removing the diagonal.

Now let $T$ be a triangulation  of the unpunctured surface $(S,M)$. If
$\tau_k\in T$ is an interior arc, then
$\tau_k $ lies in precisely two triangles in $T$, hence $\tau_k$ is
the diagonal of a unique quadrilateral $Q_{\tau_k}$ in $T$. We
associate to this quadrilateral a tile $\Sbar_k$ by assigning the weight
$x_k$ to the diagonal and the weights $x_a,x_b,x_c,x_d$ to the sides
of $\Sbar_k$ in such a way that there is a homeomorphism $\Sbar_k\to
Q_{\tau_k}$ which sends the edge with weight $x_i$ to the arc labeled
$\tau_i$, $i=a,b,c,d,k$, see Figure \ref{figtile}.

\end{subsection}

\begin{subsection}{The graph $\overline{G}_{T,\zg}$}\label{sect 1.2}
Let $T$ be a triangulation of an unpunctured surface $(S,M)$ and let
$\zg$ be an arc in $(S,M)$ which is not in $T$. 
Choose an orientation on $\zg$ and let $s\in M$ be its starting point, and let
$t\in M$ be its endpoint.
We denote by
\[ p_0=s, p_1, p_2, \ldots, p_{d+1}=t
\]
the points of intersection of $\zg$ and $T$ in order. Let $i_1,
i_2, \ldots, i_d$ be such that $p_k$ lies on the arc $\tau_{i_k}\in T$.
Note that $i_k$ may be equal to $i_j$ even if $k\ne j$. 
Let $\tilde S_1,\tilde S_2,\ldots,\tilde S_d $ be a sequence of tiles
so that $\tilde S_k$ is isomorphic to the tile $\Sbar_{i_k}$, for
$k=1,2,\ldots,d$.

For $k$ from $0$ to $d$, let $\zg_k$  denote the segment of the path
$\zg$ from the point $p_k $ to the point $p_{k+1}$. Each $\zg_k$ lies
in exactly one triangle $\zD_k$ in $T$, and if $1\le k\le d-1$ then $\zD_k$ is formed by the arcs $\tau_{i_k},
\tau_{i_{k+1}}$, and a third arc that we denote by $\tau_{[\zg_k]}$.

We will define a graph $\overline{G}_{T,\zg}$ by recursive glueing of tiles. Start
with $\overline{G}_{T,\zg,1}\cong \tilde S_1$, where we orient the tile $\tilde
S_1$ so that the diagonal goes from northwest to southeast, and the
starting point $p_0$ of $\zg$ is in the southwest corner of $\tilde
S_1$. For all $k=1,2,\ldots,d-1$ let
$\overline{G}_{T,\zg,k+1}$ be the graph obtained by adjoining the tile $\tilde
S_{k+1} $ to the tile $\tilde S_k$ of the graph $\overline{G}_{T,\zg,k}$ along
the edge weighted $x_{[\zg_k]}$, see Figure \ref{figglue}. We always
orient the tiles so that the diagonals go from northwest to southeast.
Note that the edge weighted  $x_{[\zg_k]}$ is either the northern or
the eastern edge of the tile $\tilde S_k$.
\begin{figure}
\input{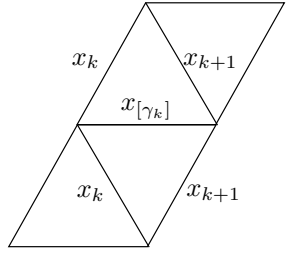}
\caption{Glueing tiles $S_k$ and $S_{k+1}$ along the edge weighted  $x_{[\zg_k]}$}\label{figglue}
\end{figure}   
Finally, we define $\overline{G}_{T,\zg}$ to be $\overline{G}_{T,\zg,d}$.

Let $G_{T,\zg}$ be the graph obtained from $\overline{G}_{T,\zg}$ by
removing the diagonal in  each tile, that is, $G_{T,\zg} $
is constructed in the same way as $\overline{G}_{T,\zg}$ but using tiles
$S_{i_k}$ instead of $\Sbar_{i_k}$.

A \emph{perfect matching} of a graph is a subset of the edges so that
each vertex is covered exactly once.  We define the weight $w(M)$ of a
perfect matching $M$ to be the product of the weights of all 
edges in $M$. 

\end{subsection}

\begin{subsection}{Cluster expansion formula}\label{sect cluster
    expansion formula}
Let $(S,M)$ be an  unpunctured surface with triangulation $T$, and let $\mathcal{A}=\mathcal{A}(\mathbf{x}_T,\mathbf{y}_T,B)$ be
the cluster algebra with principal coefficients in the initial seed 
$(\mathbf{x}_T,\mathbf{y}_T,B)$
defined in section \ref{sect FST}.
Each cluster variable in $\mathcal{A}$ corresponds to an
arc in $(S,M)$. Let $x_\zg$ be an arbitrary cluster variable
corresponding to an arc $\zg$. Choose an orientation of $\zg$, and let $\tau_{i_1}$, $\tau_{i_2} \dots,
\tau_{i_d}$ be the arcs of the triangulation that are crossed by
$\zg$ in this order, with multiplicities possible. Let $G_{T,\zg}$ be
the graph constructed in section \ref{sect 1.2}.

\begin{thm}\label{thm main} With the above notation
\[x_\zg = \sum_M \frac{w(M)\,y(M)}{x_{i_1}x_{i_2}\ldots x_{i_d}}, \]
where the sum is over all perfect matchings $M$ of $G_{T,\zg}$,  $w(M)$ is the weight of $M$, and $y(M)$ is a monomial in $\mathbf{y}_T$.
\end{thm}

The proof of Theorem \ref{thm main} will be given in Section \ref{sect proof}.

\end{subsection}

\end{section}

%%%%%%%%%%%%%%%%%%%%%%%%%%%%%%%%
%% 
%%  SECTION
%%
%%%%%%%%%%%%%%%%%%%%%%%%%%%%%%%%

\begin{section}{Proof of Theorem \ref{thm main}}\label{sect proof}
Throughout this section, $T$ is a triangulation of an
unpunctured surface $(S,M)$, $\zg$ is an arc in $S$ with a fixed orientation, and  $s\in M$ is its starting point and 
$t\in M$ is its endpoint.
Moreover, $ p_0=s, p_1, p_2, \ldots, p_{d+1}=t$ are
the points of intersection of $\zg$ and $T$ in order, and $i_1,
i_2, \ldots, i_d$ are such that $p_k$ lies on the arc $\tau_{i_k}\in T$.

\begin{subsection}{Complete $(T,\zg)$-paths}
Following \cite{ST}, 
we will consider paths $\za$ in $S$ that are concatenations of arcs in the
triangulation $T$, more precisely, $\za=
(\za_1,\za_2,\ldots,\za_{\ell(\za)})$ with 
 $\za_i \in T$, for $i=1,2,\ldots, \ell(\za)$ and the starting point
of  $\za_i$ is the endpoint of $\za_{i-1}$. Such
 a path is called a \emph{$T$-path}. 

 We call a
 $T$-path $\za=(\za_1,\za_2,\dots ,
\za_{\ell(\za)}) $  a \emph{complete $(T,\gamma)$-path}
if the following axioms hold: 

\begin{itemize}
\item[(T1)] The even arcs are precisely the arcs crossed by $\zg$ in
  order, that is, $\za_{2k}=\tau_{i_k}$.
\item[(T2)] For all $k=0,1,2,\ldots,d$, the segment $\zg_k$ is
  homotopic to the segment of the path $\za$ starting at the point
  $p_k$ following $\za_{2k},\za_{2k+1}$ and $\za_{2k+2}$ until the
  point $p_{k+1}$.
\end{itemize}

We define the Laurent monomial $x(\za)$ of the complete $(T,\zg)$-path $\za$ by
\[x(\za)=\prod_{i \textup{ odd}} x_{\za_i}\prod_{i \textup{ even}} x_{\za_i}^{-1}. \]
 \begin{rem}\begin{itemize}
\item Every complete $(T,\zg)$-path starts and ends at the same point
  as $\zg$, because of \textup{(T2)}.
\item Every complete $(T,\zg)$-path has length $2d+1$.
\item  For all arcs $\tau$ in the triangulation $T$, the number of
  times that
$\tau$ occurs as $\alpha_{2k}$ is exactly the number
  of crossings between $\zg$ and $\tau$.
\item In contrast to the ordinary $(T,\zg)$-paths defined in \cite{ST}, complete $(T,\zg)$-paths allow backtracking.
\item The denominator of the Laurent monomial $x(\za)$ is equal to $x_{i_1}x_{i_2}\cdots x_{i_d}$.
\end{itemize}
\end{rem}

\end{subsection}

\begin{subsection}{Universal cover}\label{sect universal cover}
Let $\pi:\tilde S \to S$ be a universal cover of the surface $S$, and
let $\tilde M=\pi^{-1}(M)$ and $\tilde T=\pi^{-1}(T)$. 

Choose $\tilde s \in \pi^{-1}(s)$. There exists a unique lift $\tilde
\zg$ of $\zg$ starting at $\tilde s$. Then $\tilde \zg$ is the
concatenation of subpaths $\tilde \zg_0,\tilde \zg_1,\ldots,\tilde
\zg_{d+1}$ where $\tilde \zg_k$ is a path from a point $\tilde p_k$ to
a point $\tilde p_{k+1}$ such that $\tilde \zg_k$ is a lift of $\zg_k$
and $\tilde p_k\in\pi^{-1}(p_k)$, for $k=0,1,\ldots, d+1$.
Let $\tilde t=\tilde p_{d+1}\in \pi^{-1}(t)$.

For $k$ from $1$ to $d$, let $\tilde \tau_{i_k}$ be the unique lift of
$\tau_{i_k}$ running through $\tilde p_k$ and let $\tilde
\tau_{[\zg_k]}$ be the  unique lift of $\tau_{[\zg_k]}$ that is
bounding a triangle in $\tilde S$ with $\tilde \tau_{i_k}$ and $\tilde
\tau_{i_{k+1}}$.
Each $\tilde \zg_k$ lies in exactly one triangle $\tilde \zD_k$ in
$\tilde T$. 
Let $\tilde S(\zg)\subset \tilde S$ be the union of the triangles
$\tilde \zD_0, \tilde \zD_1,\ldots,\tilde \zD_{d+1}$ and let $\tilde
M(\zg)=\tilde M\cap \tilde S(\zg)$ and $\tilde T(\zg)=\tilde T\cap
\tilde S(\zg)$.
Then $(\tilde S(\zg),\tilde M(\zg))$ is a simply connected unpunctured
surface of which $\tilde T(\zg) $ is a triangulation. This
triangulation $\tilde T(\zg)$ consists of arcs $\tilde
\tau_{i_k},\,\tilde\tau_{[\zg_k]}$ with $k=1,2,\ldots,d,$ and two arcs
incident to $\tilde s$ and two arcs incident to $\tilde t$.

The underlying graph of $\tilde T(\zg)$ is the graph with vertex set
$\tilde M(\zg)$ and whose set of edges consists of the (unoriented)
arcs in $\tilde T(\zg)$.

By \cite[Section 5.5]{S2}, we can compute the Laurent expansion of $x_{\zg}$ 
using complete $(\tilde T(\zg),\tilde \zg)$-paths in $(\tilde S(\zg),\tilde M(\zg))$. 
\end{subsection} 

\begin{subsection}{Folding}\label{sect folding}
The graph $\overline{G}_{T,\zg}$ was constructed by glueing tiles $\tilde S_{k+1}$ to
tiles $\tilde S_k$ along edges with weight 
$x_{[\zg_k]}$, see figure \ref{figglue}. 
Now we will fold the graph along the edges weighted
$x_{[\zg_k]}$, thereby identifying the two triangles incident to
$x_{[\zg_k]}$, $k=1,2,\ldots,d-1$.

To be more precise, the edge with weight $x_{[\zg_k]}$, that lies in
the two tiles $\tilde S_{k+1}$ and $\tilde S_k$, is contained in
precisely two  triangles $\zD_{k}$ and $\zD_k'$ in $\overline{G}_{T,\zg}$:
$\zD_{k}$ lying inside the tile $\tilde S_{k}$ and $\zD_k'$ lying
inside the tile $\tilde S_{k+1}$. Both  $\zD_{k}$ and $\zD_k'$ have
weights $x_{[\zg_k]}$, $x_k$, $x_{k+1}$, but opposite
orientations. Cutting $\overline{G}_{T,\zg}$ along the edge with weight
$x_{[\zg_k]}$, one obtains two connected components. Let $R_k$ be the
component that contains the tile $\tilde S_k$ and $R_{k+1}$ the
component that contains $\tilde S_{k+1}$. 

The \emph{folding of the graph $\overline{G}_{T,\zg}$ along $x_{[\zg_k]}$} is the
graph obtained by flipping $R_{k+1}$ and then glueing it to $R_k$ by
identifying the two triangles $\zD_k $ and $\zD_k'$.

The graph obtained by consecutive folding of $\overline{G}_{T,\zg}$ along all
edges with weight $x_{[\zg_k]}$ for $k=1,2,\ldots,d-1$, is isomorphic
to the underlying graph of the triangulation $\tilde T(\zg)$ of the
unpunctured surface $(\tilde S(\zg),\tilde M(\zg))$. 
Indeed, there clearly is a bijection between the triangles in both graphs,
and, in both graphs the way the triangles are glued together is
uniquely determined by $\zg$.
 Note that the two graphs may have opposite orientations.

We obtain a map that we call the \emph{folding map}
\[
\begin{array}{rcccccc}
\phi & : & \left\{ 
\genfrac{}{}{0pt}{}{\textup{perfect matchings} }
{\textup{in $G_{T,\zg}$}} 
\right\} 
&\to & \left\{ \genfrac{}{}{0pt}{}{\textup{complete $(\tilde T(\zg),\tilde
    \zg)$-paths}} {\textup{in  $(\tilde S(\zg),\tilde M(\zg))$}}\right\} 
\\
\\
&& M&\mapsto &\tilde \za_M
\end{array}  \] 
as follows.
First we associate a path $\za_M$ in $\overline G_{T,\zg}$ to the matching $M$, by inserting a diagonal between any two consecutive edges of the perfect matching. More precisely, $\za_M$ is the path starting at $s$ going along the unique edge of $M$ that is incident to $s$, then going along the diagonal of the first tile $\tilde S_1$, then along the unique edge of $M$ that is incident to the endpoint of that diagonal, and so forth.

Since $M$ has cardinality $d+1$, the path $\za_M$ consists of $2d+1$ edges, thus $\za=(\za_1 ,\za_2,\ldots,\za_{2d+1})$.
Now we define 
$\tilde\za_M=(\tilde\za_1,\tilde\za_2,\ldots,\tilde\za_{2d+1})$  
 by folding the path $\za_M$.
Thus, if $M=\{\zb_1,\zb_3,\ldots,\zb_{2d-1},\zb_{2d+1}\}$, where the
edges are ordered according to $\zg$, then $\phi(M) = (\tilde
\za_1,\tilde \za_2,\ldots,\tilde \za_{2d+1})$, where $\tilde
\za_{2k+1}$ is the image of $\zb_{2k+1}$ under the folding and $\tilde
\za_{2k}=\tau_{i_k}$ is the arc crossing $\zg$ at $p_k$.
Then $\phi(M)$ satisfies the axiom (T1) by construction. Moreover, $\phi(M)$ satisfies the axiom (T2), because, for each $k=0,1,\ldots,d$, the segment of the path $\phi(M) $, which  starts at the point $p_k$ following $\tilde \za_{2k}, \tilde \za_{2k+1}$ and $\tilde\za_{2k+2}$ until the point $p_{k+1}$, is homotopic to the segment $\zg_k$, since both segments lie in the simply connected triangle $\zD_k$ formed by $\tau_{i_k},\tau_{i_{k+1}}$ and $\tau_{[\zg_k]}$. Therefore, the folding map $\phi$ is well defined.

Note that it is possible that $\tilde \za_k, \tilde \za_{k+1}$ is
backtracking, that is, $\tilde \za_k$ and $\tilde
\za_{k+1}$ run along the same arc $\tilde\tau \in \tilde T(\zg)$.

\end{subsection}

\begin{subsection}{Unfolding the surface}
Let $\za$ be a boundary arc in  $(\tilde S(\zg),\tilde M(\zg))$ that
is not adjacent to $\tilde s$ and not adjacent to $\tilde  t$. Then
there is a unique 
triangle $\Delta$ in $\tilde T(\zg)$ in which $\za$ is a side. The other
two sides of $\Delta$ are two consecutive diagonals, which we denote
by $\tilde\tau_j$ and $\tilde\tau_{j+1}$, see Figure \ref{figcomplete}.

By cutting  the underlying graph of $\tilde T(\zg)$  along $\tilde\tau_j$, we obtain
two pieces. Let $R_{j+1}$ denote the piece that contains $\za,
\tilde\tau_{j+1}$ and $t$. Similarly, cutting  $(\tilde S(\zg),\tilde
M(\zg))$  along $\tilde\tau_{j+1}$, we obtain two pieces, and we denote by
 $R_{j}$ the piece that contains $s,\tilde\tau_{j}$ and $\za$.

The \emph{graph obtained by unfolding along $\za$} is the graph
obtained by flipping $R_{j}$ and then glueing it to $R_{j+1}$ along
$\za$.
In this new graph, we label the edge of $R_j$ that had the label
$\tilde\tau_{j+1}$ by $\tilde\tau_{j+1}^b$ and  the edge of $R_{j+1}$ that had the label
$\tilde\tau_{j}$ by $\tilde\tau_{j}^b$, indicating that these edges are on the
boundary of the new graph, see Figure \ref{figcomplete}.

\begin{figure}
\scalebox{0.8}{\input{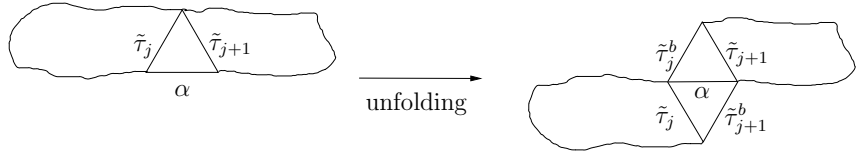}}
\caption{Completion of paths}\label{figcomplete}
\end{figure}

\begin{lem}\label{lem unfold}
The graph obtained by repeated unfolding of  the underlying graph of $\tilde T(\zg)$  along all boundary edges not adjacent to $s$ or $t$ is
isomorphic to the
graph $\overline{G}_{T,\zg}$. Moreover, for each unfolding along an edge $\za$, the
edges labeled $\tilde\tau_j^b, \tilde\tau_{j+1}^b$ are on the boundary
of $\overline{G}_{T,\zg}$ and carry the weights $x_j, x_{j+1}$ respectively, 
the edges labeled $\tilde\tau_j,\tilde\tau_{j+1}$ are diagonals in
$\overline{G}_{T,\zg}$  and carry the weights $x_j, x_{j+1}$ respectively, and
$\za$ is an interior edge of $\overline{G}_\zg$ that is not a diagonal  and carries the weight $x_{[\zg_j]}$.
\end{lem}  
\begin{pf}
This follows from the construction. \qed
\end{pf}  

\end{subsection}

\begin{subsection}{Unfolding map}
We define a map
\[
\begin{array}{rcl}
\{\textup{complete } (\tilde T(\zg),\tilde \zg)-\textup{paths}\} &\to &
\{\textup{perfect matchings of } G_{T,\zg}\}\\
\tilde\za=(\tilde\za_1,\tilde\za_2,\ldots,\tilde\za_{2d+1}) &\mapsto &
M_{\tilde \za}
=\{\zb_1,\zb_3,\zb_5,\ldots,\zb_{2d+1}\} 
\end{array}  \]

where $\zb_1=\tilde\za_1$, $\zb_{2d+1}=\tilde\za_{2d+1}$ and 
\[ \zb_{2k+1}=\left\{
\begin{array}{ll}
\tilde\za_{2k+1}&\textup{if $\tilde\za_{2k+1}$ is a boundary arc in $\tilde T(\zg)$,}\\
\tilde\tau_j^b &\textup{if $\tilde\za_{2k+1}=\tilde\tau_j$ is a diagonal in $\tilde T(\zg)$.}
\end{array}  \right.
\]

We will show that this map is well-defined. Suppose $\zb_{2k+1}$ and
$\zb_{2\ell+1}$ have a common endpoint $x$. Then  $\za_{2k+1}$ and
$\za_{2\ell+1}$ have a common endpoint $y$ in $(\tilde S(\zg),\tilde
M(\zg))$ and the two edges are not separated in the unfolding
described in Lemma \ref{lem unfold}. Consequently, there is no
triangle in $\tilde T(\zg)$ that is contained in the subpolygon spanned by
$\za_{2k+1}$ and $\za_{2\ell+1}$, hence $\za_{2k+1}$ is equal to
$\za_{2l+1}$. This implies that every arc in the subpath
$(\za_{2k+1},\za_{2k+2}\ldots \za_{2\ell+1})$ is equal to the same diagonal
$\tilde\tau_j$, and the only way this can happen is when $\ell=k+1$ and
$(\za_{2k+1},\za_{2k+2}\ldots
\za_{2\ell+1})=(\tilde\tau_j,\tilde\tau_j,\tilde\tau_j)$ and both
endpoints of $\tilde \tau_j$ are incident to an interior arc other
than $\tilde \tau_j$.
In this case, $\tilde\tau_j$ bounds the two triangles
$\tilde\tau_{j-1},\tilde\tau_{j},\tilde\tau_{[\zg_{j-1}]}$ and
$\tilde\tau_{j},\tilde\tau_{j+1},\tilde\tau_{[\zg_{j}]}$ in
$\tilde T(\zg)$. Unfolding along 
$\tilde\tau_{[\zg_{j-1}]}$ and
$\tilde\tau_{[\zg_{j}]}$ will produce edges 
$\zb_{2k+1}$ and $\zb_{2\ell+1}$ that are not adjacent, see Figure
\ref{fig 2}. 
\begin{figure}
\scalebox{1}{\input{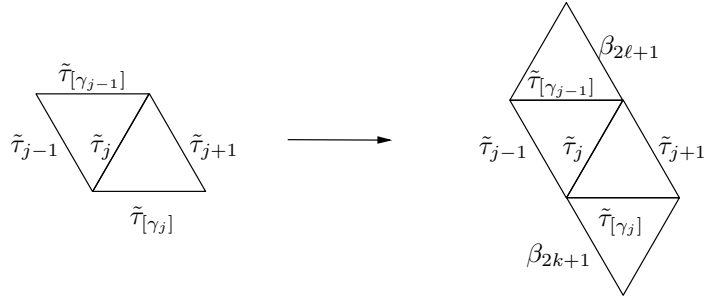}}
\caption{Unfolding along $\tilde\tau_{[\zg_{j-1}]}$ and $\tilde\tau_{[\zg_{j}]}$ }\label{fig 2}
\end{figure}    
This shows that no vertex of ${G}_{T,\zg}$ is covered twice in $M_\za$.

To show that every vertex of $G_{T,\zg}$ is covered in $M_\za$, we use a
counting argument. Indeed, the number of vertices of $G_{T,\zg}$ is
$2(d+1)$, and, on the other hand,  $2d+1$ is the length of $\za$,
since $\za$ is 
complete, and thus $M_{\tilde \za}$ has $d+1$ edges. The statement 
follows since every $\zb_j\in M_{\tilde\za}$ has two distinct endpoints.
This shows that $M_{\tilde\za} $ is a perfect matching and our map is well-defined.

\begin{lem}\label{lem bijections} The  unfolding map $\tilde\za\mapsto M_{\tilde\za}$
   is the inverse
  of the folding map $M\mapsto \tilde\za_M$. In particular, both maps are bijections.
\end{lem}  

\begin{pf} Let $\tilde \za=(\tilde \za_1,\tilde \za_2,\ldots,\tilde
  \za_{2d+1})$ be a complete $(\tilde T(\zg),\tilde \zg)$-path. Then 
$\tilde \za_{M_{\tilde \za}}=(\za_1,\za_2,\ldots,\za_{2d+1})$ where 
$\za_{2k+1}$ is the image  under folding of the arc $\tilde\tau_j^b$
  if $\tilde\za_{2k+1}=\tilde\tau_j$ is a diagonal in $\tilde T(\zg)$ or, otherwise, the image under the folding of the arc   $\tilde
  \za_{2k+1}$. Thus $\za_{2k+1}=\tilde \za_{2k+1}$. Moreover,
  $\za_{2k}=\tau_{i_k}=\tilde\za_{2k}$, and thus 
$\tilde \za_{M_{\tilde \za}}=\tilde \za$.

Conversely, let $M=\{\zb_1,\zb_3,\ldots,\zb_{2d-1},\zb_{2d+1}\} $ be a
perfect matching of $G_{T,\zg}$. Then $M_{\tilde
  \za_M}=
\{\tilde\zb_1,\tilde\zb_3,\ldots,\tilde\zb_{2d-1},\tilde\zb_{2d+1}\} $
where 
\[\begin{array}{rcl}
\tilde\zb_{2k+1}&=&\left\{
\begin{array}{rl}
\tilde\za_{2k+1} &\textup{if $\tilde\za_{2k+1} $ is a boundary arc,}\\ \\
\tilde\tau^b_{j} &\textup{if $\tilde\za_{2k+1}=\tilde\tau_j $ is a diagonal}
\end{array}  \right.  \\ 
\\
&=&\left\{
\begin{array}{rl}
\tilde\tau_{[\zg_j]} &\textup{if $\zb_{2k+1}=\tilde\tau_{[\zg_j]} $,}\\ \\
\tilde\tau^b_{j} &\textup{if $\zb_{2k+1}=\tilde\tau_j^b .$}
\end{array}  \right.\end{array}\]
Hence $M_{\tilde  \za_M}=M$.
\qed
\end{pf}  

Combining Lemma \ref{lem bijections} with the results of \cite{S2}, we obtain the following Theorem.

\begin{thm}\label{thm bijections}
There is a bijection between the set of perfect matchings of the graph $G_{T,\zg}$ and the set of complete $(T,\zg)$-paths in $(S,M)$ given by $M\mapsto \pi(\tilde \za_M)$, where $\tilde \za_M$ is the image of $M$ under the folding map and $\pi$ is induced by the universal cover $\pi:\tilde S\to S$.
Moreover, the numerator of the Laurent monomial $x(\pi(\tilde \za_M))$ of the complete $(T,\zg)$-path $\pi(\tilde \za_M)$ is equal to the weight $ w(M)$ of the matching $M$.
\end{thm}  
 
\begin{pf} The map in the Theorem is a bijection, because it is the composition of the folding map, which is a  bijection, by Lemma \ref{lem bijections}, and the map $\pi$, which is a bijection, by \cite[Lemma 5.8]{S2}. The last statement of the Theorem follows from the construction of the graph $G_{T,\zg}$.
\qed
\end{pf}

\end{subsection} 

\begin{subsection}{Proof of Theorem \ref{thm main}}\label{ssect proof}
 It has been shown in \cite[Theorem 3.2]{S2} that 

\begin{equation}\label{eq 91}
x_\zg = \sum_{\za} x(\za)\,y(\za),
\end{equation}  
where the sum is over all complete $(T,\zg)$-paths $\za$ in $(S,M)$, 
$y(\za)$ is a monomial in $\mathbf{y}_T$, and 

\begin{equation}\label{eq 92}
x(\za)=\prod_{k \ \textup{odd}} x_{\za_k} \prod_{k \ \textup{even}} x_{\za_k}^{-1}.
\end{equation}   

Applying Theorem \ref{thm bijections} to equation (\ref{eq 91}) yields

\begin{equation}\label{eq 93}
x_\zg = \sum_M w(M)\,y(M) (x_{i_1}x_{i_2}\cdots x_{i_d})^{-1},
\end{equation}  
where the sum is over all perfect matchings $M$ of $G_{T,\zg}$, $w(M)$ is the weight of the matching and $y(M)=y(\pi(\tilde \za_M))$. This completes the proof of Theorem \ref{thm main}.
\end{subsection}

\end{section}

\begin{section}{A formula for $y(M)$}\label{sect y}
In this section, we give a description of the coefficients $y(M)$ in terms of the matching $M$. First, we need to recall some results from \cite{S2}.

Recall that $T$ is a triangulation of the unpunctured surface $(S,M)$, $\zg$ is an arc in $(S,M)$ that crosses $T$ exactly $d$ times, we have
fixed an orientation for $\zg$ and denote  by  $s=p_0,p_1,\ldots,p_d,p_{d+1}=t$ 
the intersection  points of $\zg$ and $T$ in order of occurrence on
$\zg$. 
 Let $i_1,i_2,\ldots,i_d$ be such that $p_k$ lies on the arc
$\tau_{i_k}\in T$, for $k=1,2,\ldots,d$. 

For $k=0,1,\ldots,d$, let $\zg_k$ denote the segment of the path $\zg$
from  the point $p_k$ to the point $p_{k+1}$. Each $\zg_k$ lies in
exactly one triangle $\zD_k$ in $T$. If $1\le k\le d-1$, the triangle
$\zD_k$ is formed by the arcs $\tau_{i_k},
\tau_{i_{k+1}}$ and a third arc that we denote by
$\tau_{[\zg_k]}$. 

The orientation of the  surface $S$ induces an orientation on each of these triangles in such a way that, whenever two triangles $\zD,\zD'$ share an
edge $\tau$, then the orientation of $\tau $ in $\zD$ is opposite to
the orientation of $\tau $ in $\zD'$, 
There are precisely two such orientations, we assume without loss of generality that we have the ``clockwise
orientation'', that is, in each triangle $\zD$, going around the
boundary of $\zD$ according to the orientation of $\zD$ is clockwise
when looking at it from outside the surface.

Let $\za $ be a complete $(T,\zg)$-path. Then
$\za_{2k}=\tau_{i_k}$ is a common edge of the two triangles
$\zD_{k-1}$ and $\zD_k$.
We say that $\za_{2k}$ is $\zg$-\emph{oriented} if the orientation of
$\za_{2k}$ in the path $\za$ is the same as the orientation of
$\tau_{i_k}$  in the triangle $\zD_k$, see Figure \ref{fig
  gammaoriented}.
\begin{figure}
\scalebox{0.9}{\input{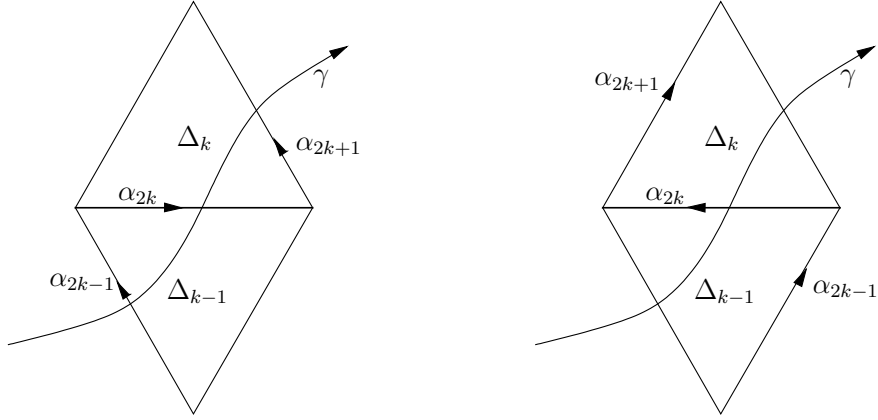}}
\caption{Two examples of the  $(T,\zg)$-path segment
  $(\za_{2k-1},\za_{2k},\za_{2k+1})$. On the left, 
  $\za_{2k}$ is not $\zg$-oriented and on the right, $\za_{2k}$ is $\zg$-oriented.}\label{fig gammaoriented}
\end{figure}   

It is shown in \cite[Theorem 3.2]{S2} that 
\begin{equation}\label{eq 28}
 y(\za)= \prod_{k:
  \za_{2k}\textup{ is $\zg$-oriented}} y_{i_k}.
\end{equation}

Each perfect matching $M$ of  $G_{T,\zg}$  induces a path $\za_M$ in $\overline{G}_{T,\zg}$ as in the construction of the folding map in section \ref{sect folding}. 
The even arcs of $\za_M$ are the diagonals of the graph $\overline{G}_{T,\zg}$. We say that an even arc of $\za_M$ has upward orientation if $\za_M$ is directed from southeast to northwest on that even arc, otherwise we say that the arc has downward orientation. 
If going upward on the first even arc of $\za_M$ is $\zg$-oriented then  
we have that
the $(2k)$-th arc  of $\pi(\tilde{\za}_M)$ is $\zg$-oriented 
if and only if the $2k$-th arc of $\za_M$ is upward if $k$ is odd, and downward if $k$ is even.
If, on the other hand,  going downward on the first even arc of $\za_M$ is 
$\zg$-oriented then  we have that
the $(2k)$-th arc  of $\pi(\tilde{\za}_M)$ is $\zg$-oriented
if and only if the $2k$-th arc of $\za_M$ is downward if $k$ is odd, and upward
 if $k$ is even.

There are precisely two perfect matchings $M_+$ and $M_-$ of $G_{T,\zg}$ that contain only boundary edges of  $G_{T,\zg}$. The orientations of the even arcs in both of the induced $(T,\zg)$-paths $\za_{M_+}$ and $\za_{M_-}$ are alternatingly upward and downward, thus for one of the two paths, say $M_+$, each even arc of $\pi(\tilde\za_{M_+})$ is $\zg$-oriented, whereas for $M_-$ none of the even arcs of $\pi(\tilde\za_{M_-})$ is $\zg$-oriented.
That is, $y(M_-)=1$ and $y(M_+)=y_{i_1}y_{i_2}\cdots y_{i_d}$.

For an arbitrary perfect matching $M$, the coefficient $y(M)$ is determined by the set of edges of the symmetric difference  
$M_-\ominus M =(M_-\cup M)\setminus (M_-\cap M)$ as follows.

\begin{thm}\label{thm y}
The set $M_-\ominus M$ is the set of boundary edges of a (possibly disconnected) subgraph $G_M$ of $G_{T,\zg}$ which is a union of tiles
\[ G_M =\cup_{j\in J} S_j. \]
Moreover,
\[y(M)=\prod_{j\in J} y_{i_j}\]
\end{thm}  

\begin{pf} Choose
any edge $e_1$ and either endpoint in $M_-\setminus(M_-\cap M)$, and walk along that edge until its other endpoint.  Since $M$ is a perfect matching, this endpoint is incident to an edge $e_2$ in in $M$, which is different  from $e_1$ and, hence, not in $M_-$. Thus $e_2\in M\setminus(M_-\cap M)$. 
Now walk along $e_2$ until its other endpoint. This endpoint is incident to an edge $e_3$ in $M_-$ which is different from $e_2$, and, hence, not in $M$. Thus $e_3\in M_-\setminus(M_-\cap M)$. Continuing this way, we construct a sequence of edges in $M_-\ominus M$. Since $G_{T,\zg}$ has only finitely many edges, this sequence  must become periodic after a certain number of steps; thus there exist $p,N$ such that 
$e_k=e_{k+p}$ for all $k\ge N$.

We will show that one can take $N=1$. Suppose to the contrary that $N\ge 2$ is the smallest integer such that $e_k=e_{k+p}$ for all $k\ge N$. Then $e_{N-1},e_N$ and $e_{N+p-1}$ share a common endpoint. But  $e_{N-1},e_N$
 and $e_{N+p-1}$ are elements of the union of two perfect matchings, hence $e_{N-1}=e_{N+p-1}$, contradicting the minimality of $N$.

Therefore the sequence $e_1, e_2,\ldots, e_p$ in $M\ominus M_-$ is the set of boundary edges of a connected subgraph of $G_{T,\zg}$ which is a union of tiles.

The graph $G_M$ is the union of these connected subgraphs and, hence, it is a union of tiles.
Let $H$ be a connected component of $G_M$. There are precisely two perfect matchings $M_-(H)$ and $M_+(H)$ of $H$ that consist only of boundary edges of $H$. Clearly, these two matchings are $M_-\cap E(H)$ and $M\cap E(H)$, where $E(H)$ is the set of edges of the graph $H$. 
Therefore, in each tile of $H$, the orientation of the diagonal in $\za_{M_-}$ and $\za_M$ are  opposite. 
The restrictions of $M_-$ and $M$ to $E(G_{T,\zg})\setminus E(G_M)$ are  identical, hence in each tile of $G_{T,\zg} \setminus G_M$, the orientations of the diagonal in $\za_{M_-}$ and $\za_M$ are equal.
 It follows from equation (\ref{eq 28}) that $y(M)=\prod_{j\in J} y_{i_j}$.
\qed
\end{pf}

It has been shown in \cite{FZ4} that, for any cluster variable $x_\zg$ in  $\mathcal{A}$, its Laurent expansion in the initial seed $(\mathbf{x}_T,\mathbf{y}_T,B_T)$ is homogeneous
 with respect to the grading given by
$\textup{deg}(x_i)=\mathbf{e}_i$ and
$\textup{deg}(y_i)=B_T\mathbf{e}_i$, where
$\mathbf{e}_i=(0,\ldots,0,1,0,\ldots,0)^T 
\in\mathbb{Z}^n$ with $1$ at
position $i$. By
definition, the \emph{$g$-vector} $g_\zg$ of a cluster variable $x_\zg$ is the 
degree of its Laurent expansion  with respect to this grading.

\begin{cor}\label{cor g}
The $g$-vector $g_\zg$ of $x_\zg$ is given by
\[g_\zg = \deg \frac{w(M_-)}{x_{i_1}x_{i_2}\cdots x_{i_d}}.\]
\end{cor}
\begin{pf} This follows from the fact that $y(M_-)=1$.
\qed
\end{pf}

%%
%% ADDED IN DESCRIPTION OF CONNECTION TO HEIGHT FUNCTIONS
%% Gregg 10/15 Turned this into a remark and slightly rewrote
%%

\begin{rem} \label{height}
The formula for $y(M)$ can also be phrased in terms of height functions.  As described in section 3 of \cite{ProppLattice}, one way to define the height function on the faces of a bipartite planar graph $G$, covered by a perfect matching $M$, is to superimpose each matching with the fixed matching $M_{\hat{0}}$ (the unique matching of minimal height).  In the case where $G$ is a snake graph, we take $M_{\hat{0}}$ to be $M_-$, one of the two matchings of $G$ only involving edges on the boundary.  Color the vertices of $G$ black and white so that no two adjacents vertices have the same color.  In this superposition, we orient edges of $M$ from black to white, and edges of $M_-$ from white to black.  We thereby obtain a spanning set of cycles, and removing the cycles of length two exactly corresponds to taking the symmetric difference $M \ominus M_-$.  We can read the resulting graph as a relief-map, in which the altitude changes by $+1$ or $-1$ as one crosses over a contour line, according to whether the counter-line is directed clockwise or counter-clockwise.  By this procedure, we obtain a height function $h_M : F(G) \rightarrow \mathbb{Z}$ which assigns integers to the faces of graph $G$.  When $G$ is a snake graph, the set of faces $F(G)$ is simply the set of tiles $\{S_j\}$ of $G$.  Comparing with the defintion of $y(M)$ in Theorem \ref{thm y}, we see that $$y(M) = \prod_{S_j \in F(G)} y_j^{h_M(j)}.$$  An alternative defintion of height functions comes from \cite{EKLP} by translating the matching problem into a domino tiling problem on a region colored as a checkerboard.  We imagine an ant starting at an arbitary vertex at height $0$, walking along the boundary of each domino, and changing its height by $+1$ or $-1$ as it traverses the boundary of a black or white square, respectively.  The values of the height function under these two formulations agree up to scaling by four.
\end{rem}

\end{section} 
\begin{section}{Cluster expansion without matchings}\label{sect F-polynomial}
In this section, we give a formula for the cluster expansion of $x_\zg$ in terms of the graph $G_{T,\zg}$ only. 

For any graph $H$, let $c(H)$ be the number of connected components of $H$. 
Let $E(H)$ be the set of edges of $H$, and denote by $\partial H$ the set of boundary edges of $H$.
Define $\mathcal{H}_k$ to be 
the set of all subgraphs $H$ of $G_{T,\zg}$ such that $H$ is a union of $k$ tiles $H=S_{j_1}\cup\cdots\cup S_{j_k}$ and 
the number of edges of $M_-$ that are contained in $H$ is equal to $k+c(H)$.
For $H\in \mathcal{H}_k$, let \[y(H)=\prod_{S_{i_j} \textup{\,tile\,in\,}H} y_{i_j}.\]

\begin{thm}\label{thm F}
The cluster expansion of the cluster variable $x_\zg$ is given by
\[x_\zg=\sum_{k=0}^d \  \sum_{H\in \mathcal{H}_k} \frac{w(\partial H\ominus M_-)\,y(H)}{x_{i_1} x_{i_2}\cdots x_{i_d}},\]
\end{thm}

\begin{pf}
It follows from the theorems \ref{thm main} and \ref{thm y} that 
\[x_\zg=\sum_{k=1}^d\quad \sum_{M:\mid y(M)\mid =k}\quad  \frac{w(M)\,y(G_M)}{x_{i_1} x_{i_2}\cdots x_{i_d}},\]
where  $| y(M)| $ is the number of tiles in $G_M$. 
We will  show that for all $k$, the map $M\mapsto G_M$ is a bijection between the set of perfect matchings $M$ of $G_{T,\zg}$ such that $| y(M)| =k$ and  the set $\mathcal{H}_k$.
\begin{itemize}
\item[-] The map is well-defined. 
Clearly, $G_M$ is the union of $k$ tiles. Moreover, $E(G_M)\cap M_-$ is a perfect matching of $G_M$, since $M_-$ consists of every other boundary edge of $G_{T,\zg}$.
Thus the cardinality of $(E(G_M)\cap M_-)$ is half the number of vertices of $G_M$, which is equal to $2k+2c(G_M)$. Therefore, the cardinality of $(E(G_M)\cap M_-)$ is $k+c(G_M)$ and $G_M\in \mathcal{H}_k$.
\item[-] The map is injective, since two graphs $G_M,G_{M'}$ are equal if and only if their boundaries are.
\item[-] The map is surjective. 
Let $H=S_{j_1}\cup\cdots\cup S_{j_k}$ such that the cardinality of
 $E(H)\cap M_-$ equals $k+c(H)$. 
The boundary of $H$ consists of $2k+2+2c(H)$ edges, 
half of which lie in $M_-$. 
As in the proof of Theorem \ref{thm y}, let $M_-(H)=E(H)\cap M_-$ and $M_+(H)$ 
be the two perfect matchings of $H$ that consist of boundary edges only.
Let $M = M_+(H)\cup (M_-\setminus M_-(H))$.
Then $M$ is a perfect matching of $G_{T,\zg}$ such that $G_M=H$, and moreover, $|y(M)|$ is equal to the number of tiles in $H$, which is $k$.
Thus the map is surjective.
\end{itemize}

Now the boundary edges of $G_M$ are precisely the elements of $M\ominus M_-$, which implies that 
$\partial(G_M)\ominus M_- = (M\ominus M_-)\ominus M_- = M\ominus (M_-\ominus M_-) =M$.
Therefore $w(M)=w(\partial(G_M)\ominus M_-)$, and this completes the proof.
\qed
\end{pf}

\begin{cor} The $F$-polynomial of $\zg$ is given by
\[F_\zg=\sum_{k=0}^d \  \sum_{H\in\mathcal{H}_k} y(H).\]
\end{cor}
\end{section}

\begin{section}{Example}\label{sect example}

\begin{figure}[h!]
\begin{center}
\input{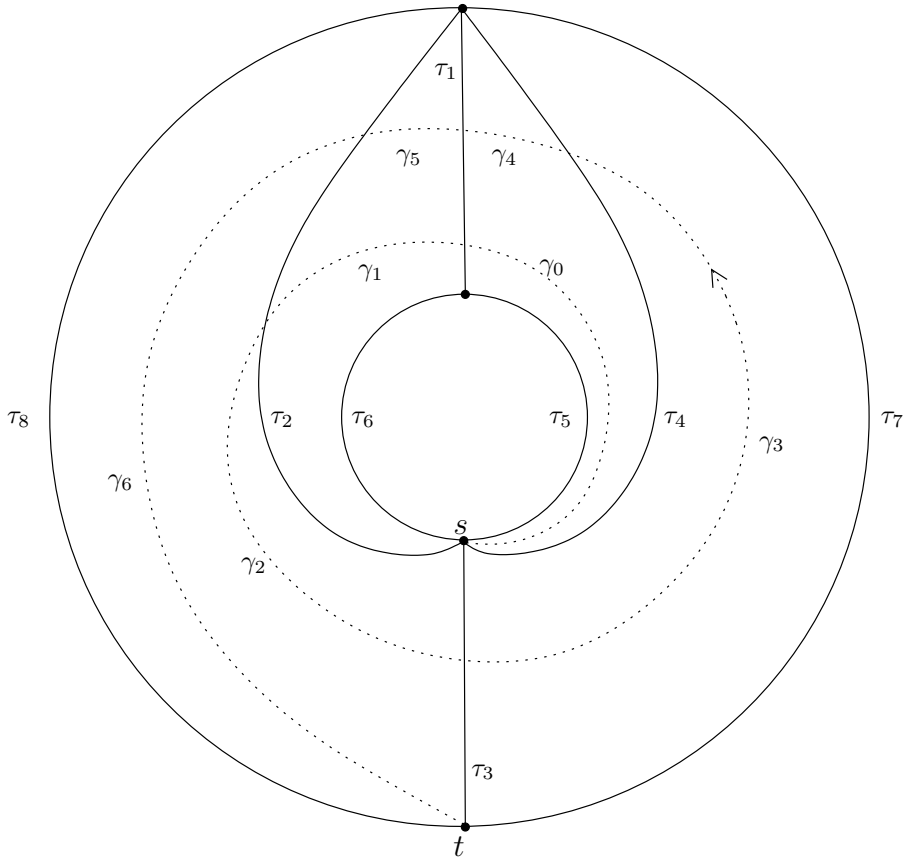}
\caption{Triangulated surface with dotted arc $\zg$} \label{figex1}
\end{center} 
\end{figure}

\begin{figure}[h!]
\begin{center}
\scalebox{0.9}{\input{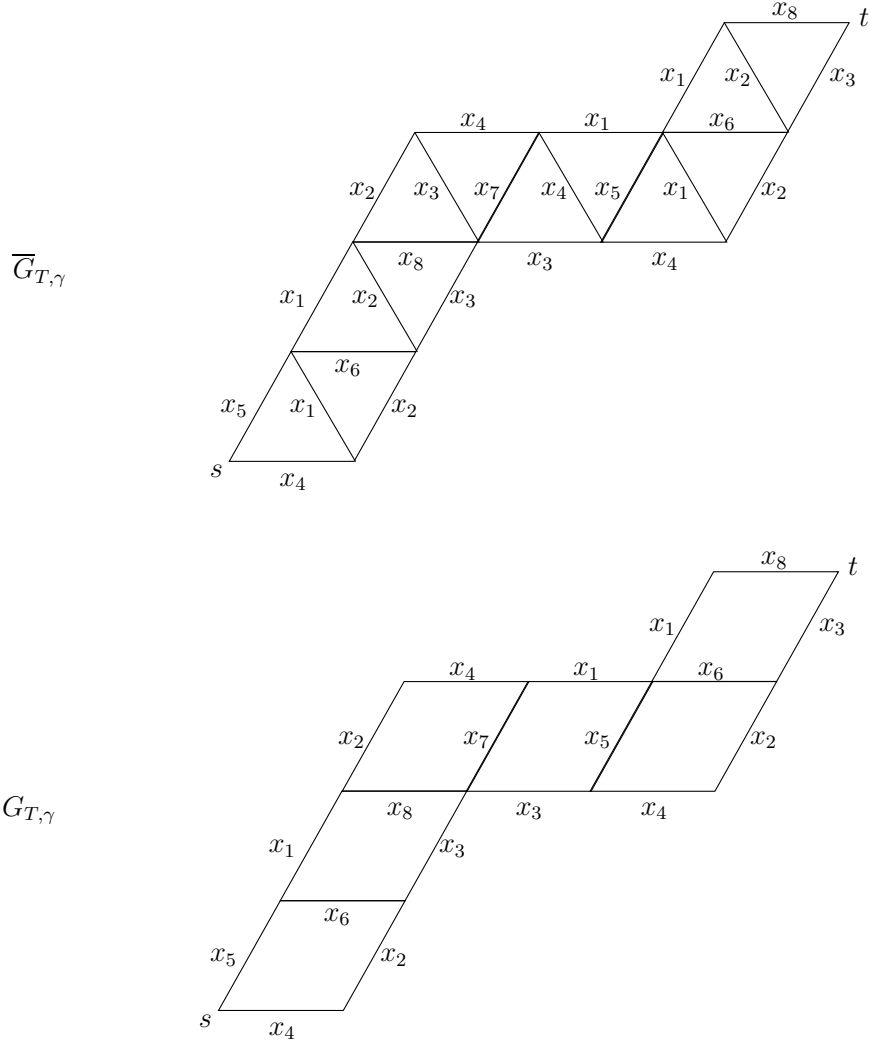}}
\caption{Construction of the graphs $\overline{G}_{T,\zg}$ and $G_{T,\zg}$ } \label{figex1a}
\end{center} 
\end{figure}

We illustrate Theorem \ref{thm main}, Theorem \ref{thm y} and Theorem \ref{thm F} in an
example. Let $(S,M)$ be the annulus with two marked points on each of
the two boundary components, and let $T=\{\tau_1,\ldots,\tau_8\}$ be
the triangulation shown in Figure \ref{figex1}. Let $\zg$ be the
dotted arc in that figure. It has $d=6$ crossings with the
triangulation. The sequence of crossed arcs 
$\tau_{i_1},\ldots,\tau_{i_6}$ is
$\tau_1,\tau_2,\tau_3,\tau_4,\tau_1,\tau_2$, and the corresponding segments
$\zg_0,\ldots,\zg_6$ of the arc $\zg$ are labeled in the figure.
Moreover, $\tau_{[\zg_1]}=\tau_6$, $\tau_{[\zg_2]}=\tau_8$,
$\tau_{[\zg_3]}=\tau_7$,  $\tau_{[\zg_4]}=\tau_5$ and
$\tau_{[\zg_5]}=\tau_6$.

The graph  $G_{T,\zg}$ is obtained by glueing the corresponding six
tiles $\tilde S_1$,  $\tilde S_2$,  $\tilde S_3$,  $\tilde S_4$,
$\tilde S_1$, and $\tilde S_2$. The result is shown in Figure \ref{figex1a}. 

Theorems \ref{thm main} and \ref{thm y} imply that $x_\zg (x_{i_1}x_{i_2}\cdots x_{i_d})$ is equal to
\[\begin{array}{cccccc}
&x_4x_6x_8x_4x_4x_6x_8 \,y_1y_3y_4y_1&+& x_4x_6x_8x_4x_4x_1x_3 \, y_1y_3y_4y_1y_2\\+&
x_4x_6x_8x_4x_5x_2x_8 \,y_1y_3y_4
&+&x_4x_6x_2x_3x_1x_2x_8 \,y_1
\\
+& x_4x_6x_2x_7x_5x_2x_8\,y_1y_4 &+& x_4x_6x_2x_7x_4x_6x_8\,y_1y_4y_1\\
+& x_4x_1x_3x_4x_4x_6x_8 \,y_1y_2y_3y_4y_1 &+& x_4x_1x_3x_4x_4x_1x_3\,y_1y_2y_3y_4y_1y_2
\\
+& x_4x_1x_3x_4x_5x_2x_8\,y_1y_2y_3y_4
&+& x_5x_2x_8x_4x_4x_6x_8\, y_3y_4y_1\\
+& x_5x_2x_8x_4x_4x_1x_3\, y_3y_4y_1y_2 &+&x_5x_2x_8x_4x_5x_2x_8\, y_3y_4\\
+& x_5x_2x_2x_3x_1x_2x_8\, &+& x_5x_2x_2x_7x_4x_6x_8 \, y_4y_1 \\
+& x_5x_2x_2x_7x_4x_1x_3\, y_4y_1y_2
&+& x_5x_2x_2x_7x_5x_2x_8\,y_4,
\end{array}\]which is equal to 

\[\begin{array}{cccccc}
 &x_4^3\,y_1^2y_3y_4&+& x_1x_3x_4^3 \, y_1^2y_2y_3y_4 \\
 +& x_2x_4^2 \,y_1y_3y_4&+&x_1x_2^2x_3x_4\,y_1\\
+& x_2^2x_4\,y_1y_4 &+& x_2x_4^2\,y_1^2y_4\\
+& x_1x_3x_4^3\,y_1^2y_2y_3y_4 &+&x_1^2x_3^2x_4^3\,y_1^2y_2^2y_3y_4\\
+& x_1x_2x_3x_4^2\,y_1y_2y_3y_4
&+& x_2x_4^2 \,y_1 y_3y_4\\
+& x_1x_2x_3x_4^2\, y_1y_2y_3y_4
&+&x_2^2x_4\, y_3y_4 \\
+& x_1x_2^3x_3 &+&x_2^2x_4 \, y_1y_4\\
+& x_1x_2^2x_3x_4\, y_1y_2y_4
&+& x_2^3\,y_4.\end{array}\]
For example, the first term corresponds to the matching $M$ consisting of the horizontal edges  of the first three tiles and the horizontal edges  of  the last two tiles. The matching $M_-$ consists in the boundary edges weighted $x_5$ and $x_2$ in the first tile, $x_2$ in the third tile, $x_1$ and $x_3 $ in the forth, $x_2$ in the fifth and $x_8$ in the sixth tile.
Thus $M_-\ominus M =(M_-\cup M)\setminus (M_-\cap M)$ is the union of the first, third, forth and fifth tile, whence $y(M)=y_{i_1}y_{i_3}y_{i_4}y_{i_5}=y_1y_3y_4y_1$.

To illustrate Theorem \ref{thm F}, let $k=2$. Then $\mathcal{H}_k$ consists of the subgraphs $H$ of $G_{T,\zg}$ which are unions of two tiles and such that $E(H)\cap M_-$ has three elements if $H$ is connected, respectively four elements if $H$ has two connected components.
Thus $\mathcal{H}_2$ has three elements
\[
\mathcal{H}_2=\{S_{i_3}\cup S_{i_4}, S_{i_4}\cup S_{i_5}, S_{i_1}\cup S_{i_4}\}
\]
corresponding to the three terms
\[x_2^2x_4y_3y_4 , x_2^2x_4y_1y_4 \textup{ and }x_2^2x_4y_1y_4. \]
\end{section} 

{}

\vspace{1cm}
\begin{tabular}{ll}{\small
Gregg Musiker }& {\small Ralf Schiffler}\\
{\small Department of Mathematics, Room 2-332 }\qquad\qquad& {\small  Department of Mathematics }\\
{\small Massachusetts Institute of Technology} &{\small University of Connecticut}\\
{\small 77 Massachusetts Ave.}&{\small 196 Auditorium Road}\\
{\small Cambridge, MA 02139} &{\small  Storrs, CT 06269-3009} \\
{\small musiker (at) math.mit.edu }&{\small schiffler (at) math.uconn.edu}
\end{tabular}

\end{document}